\documentclass{article}
\usepackage[utf8]{inputenc}
\usepackage{amssymb}
\usepackage{amsmath}
\usepackage{amsthm}
\usepackage{enumerate}
\usepackage{mathrsfs}
\usepackage{microtype}
\usepackage{dsfont}
\usepackage{bm}
\usepackage{dsfont, graphicx, url}

\newtheorem{theorem}{Theorem}

\newcommand{\ignore}[1]{}

\usepackage[top=1in, bottom=1in]{geometry}

\newcommand{\myAuthor}[3]{\noindent {\scshape #1}\\ \texttt{#2}\\ #3}

\begin{document}
\title{Connected-Intersecting Families of Graphs}
\date{January 5, 2019}

\author{Aaron Berger \and Ross Berkowitz \and Pat Devlin \and Michael Doppelt \and Sonali Durham \and Tessa Murthy \and Harish Vemuri}

\maketitle
\begin{abstract}
For a graph property $\mathcal{P}$ and a common vertex set $V = \{1, 2, \ldots, n\}$, a family of graphs on $V$ is \emph{$\mathcal{P}$-intersecting} iff $G \cap H$ satisfies $\mathcal{P}$ for all $G,H$ in the family.  Addressing a question of Chung, Graham, Frankl, and Shearer, we explore---for various $\mathcal{P}$---the maximum cardinality among all $\mathcal{P}$-intersecting families of graphs.  In the connected-intersecting case, we resolve the question completely by a short linear algebraic proof showing this maximum is attained by taking all graphs containing a fixed spanning tree (though we show other extremal constructions as well).  We also present a new lower bound for containing unions of a fixed subgraph.
\end{abstract}

\renewcommand{\thefootnote}{\fnsymbol{footnote}}

\footnotetext{AMS 2010 subject classification: 05C35, 05C40, 05D05, 05C50, 90C27}
\footnotetext{Key words and phrases:  Erd\H{o}s--Ko--Rado, graph families, intersection, connectivity, anti-clusters}

\section{Introduction}
One of the earliest and most important results in extremal combinatorics is the Erd\H{o}s--Ko--Rado theorem \cite{ekr}, which states the following:
\begin{theorem}[Erd\H{o}s--Ko--Rado \cite{ekr}]
Suppose $n \geq 2k$ and $\mathcal{F}$ is a collection of $k$-element subsets of $\{1, 2, \ldots , n\}$ such that for all $A, B \in \mathcal{F}$, $A \cap B \neq \emptyset$.  Then $|\mathcal{F}| \leq {n-1 \choose k-1}$.  Moreover, for $n > 2k$ equality is attained iff $\mathcal{F}$ is a collection of all $k$-element sets containing some particular element.
\end{theorem}
\noindent Since its publication in 1961, the Erd\H{o}s--Ko--Rado theorem has been generalized in many ways (e.g., via the equally classic Hilton--Milner theorem \cite{hm}), and it has spawned a host of related questions (see \cite{permutations} as an example and \cite{survey, yuvalsurvey} for two recent surveys).  In a 1986 publication \cite{cgfs}, Chung, Graham, Frankl, and Shearer posed a broad generalization of the problem requiring pairwise intersections to have certain \textit{structure}.  Of this, we recall only the case for graphs, which is our primary focus.

Suppose $\mathcal{P}$ is a graph property and $\mathcal{F}$ is a family of graphs each on some common vertex set $V = \{1, 2, \ldots, n\}$.  We say $\mathcal{F}$ is \emph{$\mathcal{P}$-intersecting} iff $G \cap H$ satisfies $\mathcal{P}$ for all $G,H \in \mathcal{F}$.  In this language, the main question posed in \cite{cgfs} is the following:

\medskip
\noindent \textbf{Question A:} For a fixed property $\mathcal{P}$, how large is the largest $\mathcal{P}$-intersecting family of graphs?
\medskip

To establish notation, let
\[
\mu (\mathcal{P}, n) := \max \left \{|\mathcal{F}|2^{-{n \choose 2}} \ : \ \text{$\mathcal{F}$ is $\mathcal{P}$-intersecting with vertex set $\{1, 2, \ldots, n\}$}\right\}
\]
denote the maximum size of such a family.  Note this count is normalized in our definition so that $\mu(\mathcal{P},n) \in [0,1]$ is the fraction of graphs on $\{1, 2, \ldots, n\}$ in $\mathcal{F}$.  Thus, Question A is to bound $\mu(\mathcal{P}, n)$.

Although Question A was posed in full generality, the majority of research to date has been for properties of the form ``contains a copy of the graph $\Gamma$."  Graph families for which each pairwise intersection contains some copy of $\Gamma$ are called $\Gamma$-intersecting, and triangle-intersecting families were of particular interest.  For any fixed (non-empty) graph $\Gamma$, we have the almost immediate bounds
\begin{equation}\label{sillyBounds}
2 ^{-|E(\Gamma)|} \leq \mu(\text{contains $\Gamma$}, n) \leq 1/2.
\end{equation}
The upper bound holds since an intersecting family cannot contain both a graph and its complement, and the lower bound is attained by fixing a particular copy of $\Gamma$ and taking all graphs containing it.

In 1976, Simonovits and S\'os conjectured that for triangle-intersecting families, the lower bound is actually correct.  In fact, this conjecture was a motivation for \cite{cgfs}, who proved (\textit{inter alia})
\[
1/8 \leq \mu(\text{contains $\triangle$}, n) \leq \mu(\text{is not bipartite}, n) \leq 1/4.
\]
Their argument relies on entropy and is (essentially) the earliest use of what is now known as \textit{Shearer's lemma}.  This upper bound remained unchanged for decades until a breakthrough of Ellis, Filmus, and Freidgut \cite{eff} using Fourier analytic techniques to show
\[
\mu(\text{contains $\triangle$}, n) = \mu(\text{is not bipartite}, n) = 1/8.
\]
This result completely settles the case for triangles, and \cite{eff} also shows that equality is attained only for the canonical (i.e., `na\"ive') construction described above.  They also obtain an analogous statement about triangle-intersecting families under more general product measures.

For other intersecting families, it is a folklore result---alluded to in \cite{cgfs} and fully stated by Alon and Spencer \cite{as}---that if $C$ is a constellation (i.e., a forest in which each connected component is a star) then $\lim_{n \to \infty} \mu(\text{contains $C$}, n) = 1/2$.  Alon and Spencer conjecture (as echoed in \cite{eff}) that $\mu(\text{contains $\Gamma$}, n) < 0.4999$ for every other graph $\Gamma$.  This would follow from the case $\Gamma=P_3$, the path with three edges, but in fact \textit{there are no known upper bounds for any bipartite} $\Gamma$ that are asymptotically better than \eqref{sillyBounds}.  And although it is natural to conjecture that in general the lower bound of \eqref{sillyBounds} is tight, this was disproven by Christofides \cite{lowerBound} who found a construction exhibitting
\[
2^{-3} < 17/128 \leq \mu(\text{contains $P_3$}, n).
\]

\noindent For other properties, nothing else appears on the problem except for an exercise in Alon and Spencer \cite{as} showing (via a modification of the entropy proof in \cite{cgfs})
\[
\mu(\text{contains a perfect matching}, n) \leq \mu(\text{has no isolated vertices}, n) \leq 2^{-n/2},
\]
which is tight for even $n$ by taking all graphs containing some particular perfect matching.  For any fixed $r$, another modification of \cite{cgfs} also gives the upper bound $\mu(\text{not $r$-partite}, n) \leq 2^{-r}$, although we could not find any explicit reference for this.

Our main result is proven in Section \ref{section:connected}, where we consider connected-intersecting families, settling the corresponding case of Question A completely.

\begin{theorem}\label{connected}
Suppose $\mathcal{F}$ is a connected-intersecting family of graphs on $\{1, 2, \ldots, n\}$ (i.e., each pairwise intersection is connected).  Then $|\mathcal{F}| \leq 2^{{n \choose 2}} / 2^{n-1}$.  Therefore, $\mu(\text{connected}, n) = 2^{-n+1}$.
\end{theorem}

\noindent Although this is tight via the canonical examples (i.e., all graphs containing some fixed spanning tree), we also exhibit many other graph families that are just as large.  Our proof is a short linear algebraic argument, related to the idea of \textit{anti-clusters} as introduced by Griggs and Walker \cite{griggs}.

In Section \ref{section:winkler}, we also exhibit a recursive construction allowing us to find a large family of graphs for which the lower bound in \eqref{sillyBounds} is exponentially far from the truth.
\begin{theorem}\label{unions}
For any non-empty graph $\Gamma$, let $\Gamma_t = \Gamma \cup \Gamma \cup \cdots \cup \Gamma$ be the graph consisting of $t$ vertex-disjoint copies of $\Gamma$.  For all $\varepsilon > 0$, there exists an integer $T$ such that for each fixed $t > T$
\[
\lim_{n \to \infty} \mu(\text{contains $\Gamma_t$}, n) > \left((1- \varepsilon) \lim_{n \to \infty} \dfrac{\mu(\text{contains $\Gamma$}, n)}{1-\mu(\text{contains $\Gamma$}, n)}\right)^{t}.
\]
\end{theorem}
\noindent Note that for any strongly increasing property $\mathcal{P}$ (i.e., any property retained both by the addition of edges as well as the addition of vertices) the function $\mu(\mathcal{P},n)$ is non-decreasing in $n$, so the limits in Theorem \ref{unions} exist for any fixed $\Gamma$.  We conclude in Section \ref{section:conclusion} with open questions and conjectures.

\vspace*{12pt}
\noindent \textbf{Acknowledgement:} We would like to thank the support and funding of the 2018 \textit{Summer Undergraduate Math Research at Yale} (SUMRY) program, where much of this project was completed.

\section{Connected intersections}\label{section:connected}
Our proof of Theorem \ref{connected} uses anti-clusters (namely, cosets of $W$ as below) introduced in \cite{griggs}; however, instead of recalling any of this machinery, it is easier just to reproduce everything from scratch.

\begin{proof}
Suppose $\mathcal{F}$ is a connected-intersecting family.  Following a standard approach (see, e.g., \cite{cgfs}), we equip the space of graphs on $\{1, 2, \ldots, n\}$ with an operation $\oplus$ where $G \oplus H$ denotes the graph on $\{1, 2, \ldots n\}$ whose edge set is $E(G\oplus H) = E(G) \triangle E(H) = [E(G) \setminus E(H)] \cup [E(H) \setminus E(G)]$.  This gives the set of graphs on $\{1, 2, \ldots, n\}$ the structure of an ${n \choose 2}$-dimensional vector space over $\mathbb{F}_2$.

For $S \subseteq \{1, 2, \ldots, n\}$, let $B(S)$ denote the complete bipartite graph on $\{1, 2, \ldots, n\}$ with edges $u \sim v$ iff $|\{u,v\} \cap S| = 1$.  Note that the set $W = \{B(S) \ : \ S \subseteq \{1, 2, \ldots, n\} \}$ consisting of all complete bipartite graphs is a subspace: in fact $B(S) \oplus B(T) = B(S \triangle T)$.  Moreover, $|W| = 2^{n-1}$ since $B(S) = B(T)$ iff $S = T$ or $S = T^{c}$.

Finally, we claim that $\mathcal{F}$ intersects each coset of $W$ at most once.  Otherwise, $\mathcal{F}$ has distinct elements of the form $G_1 = H \oplus B(S_1)$ and $G_2 = H \oplus B(S_2)$.  Because $G_1 \neq G_2$, we know that $G_1 \oplus G_2 = B(S_1 \triangle S_2)$ is the complement of two non-empty cliques, and hence $(G_1 \oplus G_2)^c$ is disconnected.  But for any graphs $G_1 \cap G_2 \subseteq (G_1 \oplus G_2)^{c}$, which would imply $G_1 \cap G_2$ is disconnected.  Thus, since $\mathcal{F}$ intersects each coset of $W$ at most once, $|\mathcal{F}| \leq 2^{{n \choose 2}} / |W| = 2^{{n \choose 2}} / 2^{n-1}$.
\end{proof}
If $T$ is any spanning tree, then the proportion of graphs containing $T$ is precisely $2^{-|E(T)|} = 2^{-n+1}$, which yields a connected-intersecting family of maximum possible size.  However, there are many more families of this same size.  For example, let $A$ and $B$ be vertex-disjoint trees whose union spans every vertex of $\{1, 2, \ldots, n\}$ and let $S$ be a subset of edges between $A$ and $B$ such that $|S|$ is odd.  Take $\mathcal{F}$ to be the family of graphs containing $A \cup B$ as well as at least half the edges of $S$.  Then $\mathcal{F}$ is connected-intersecting since the intersection of any two elements contains $A \cup B$ as well as at least one edge of $S$ (since $|S|$ is odd), and $|\mathcal{F}|$ is exactly $2^{-|E(A \cup B)|} \cdot \frac{1}{2} = 2^{-n+1}$ of the graphs.  This construction can also be iterated to obtain large, strange families of connected-intersecting graphs.  Moreover, it is not difficult to find extremal connected-intersecting families not of this type on as few as four vertices (e.g., for $n=4$, take a $4$-cycle and also all graphs having at least 5 edges).

\section{Containing disjoint unions}\label{section:winkler}
We now turn our attention to a proof of Theorem \ref{unions}.  In keeping with the notation of the theorem, let $\Gamma$ be a fixed non-empty graph and let $\Gamma_k$ denote the graph consisting of $k$ vertex-disjoint copies of $\Gamma$.  Our construction relies on the following simple ``tensoring" trick.

Namely, suppose $\mathcal{F}$ is a graph family on $\{1, 2, \ldots, n\}$, and consider a graph $G$ on $\{1, 2, \ldots, M n\}$.  For $1 \leq i \leq M$, let $G_i$ be the subgraph of $G$ induced on the vertex set $\{v \ : \ 1 + n(i-1)\leq v \leq n + n(i-1)\}$.  As a slight abuse of notation, we say $G_i \in \mathcal{F}$ iff the graph isomorphic to $G_i$ obtained by subtracting $n(i-1)$ from the label of each vertex is an element of $\mathcal{F}$.

Finally, for integers $0 \leq a \leq b$, let $\mathcal{F}(a,b)$ denote the family of graphs $G$ on $\{1, 2, \ldots, bn\}$ such that $|\{ i \ : \ G_i \in \mathcal{F}\}| \geq a$ (i.e., graphs for which at least $a$ of the $b$ projections belong to $\mathcal{F}$).

Suppose $\mathcal{F}$ is an optimal $\Gamma$-intersecting family of graphs on $\{1, 2, \ldots, n\}$.  Then for any $x \geq 0$, $\mathcal{F}(t+x, t+2x)$ is a $\Gamma_t$-intersecting family (since for any two elements of $\mathcal{F}(t+x,t+2x)$, there are at least $t$ indices for which both projections are elements of $\mathcal{F}$).  Therefore, we have
\begin{eqnarray*}
\mu(\text{contains $\Gamma_t$}, n(t+2x)) &\geq& |\mathcal{F}(t+x, t+2x)| / 2^{{(t+2x)n \choose 2}}\\
&=& \sum_{j=t+x} ^{t+2x} {t +2x \choose j} \left( \dfrac{|\mathcal{F}|}{2^{n(n-1)/2}} \right)^j \left(1- \dfrac{|\mathcal{F}|}{2^{n(n-1)/2}} \right)^{t+2x - j}\\
&=& \sum_{j=t+x} ^{t+2x} {t +2x \choose j} \mu(\text{contains $\Gamma$}, n) ^j \left(1- \mu(\text{contains $\Gamma$}, n) \right)^{t+2x - j}.
\end{eqnarray*}
Setting $p = \lim_{n \to \infty} \mu(\text{contains $\Gamma$}, n)$ and taking limits of the above as $n \to \infty$, we obtain
\begin{equation}\label{equation:union}
\lim_{n \to \infty} \mu(\text{contains $\Gamma_t$}, n) \geq \sum_{j=t+x} ^{t+2x} {t +2x \choose j} p ^j \left(1- p \right)^{t+2x - j} \geq {t +2x \choose t+x} p ^{t+x} \left(1- p \right)^{x}.
\end{equation}
For notational ease, set $\gamma = x/(t + 2x)$ (viewing $t$ as large but arbitrary and $\gamma$ to be determined).  By a standard estimate from information theory (see e.g., \cite{informationTheory}), the right-hand side of \eqref{equation:union} is bounded by

\begin{eqnarray*}
{t +2x \choose t+x} p ^{t+x} \left(1- p \right)^{x} &=&  p^t {t +2x \choose x} p^x (1-p)^x \geq \dfrac{p^t}{t+2x+1} \left[2^{-\lg(\gamma) - (1-\gamma) \lg(1-\gamma) / \gamma} p(1-p)\right]^x\\
&=& \left[\left(\dfrac{1}{t+2x+1} \right)^{1/t} \left(\dfrac{p(1-p)}{\gamma(1-\gamma)} \right)^{\gamma /(1 - 2 \gamma)}\left(\dfrac{1-p}{1-\gamma} \right) \dfrac{p}{1-p}\right]^{t},
\end{eqnarray*}
where $\lg = \log_2$ denotes the binary logarithm.  For each value of $t$, we then pick $x$ so that $\gamma$ is as close as possible to $p \in (0, 1/2]$.  The fact that $\gamma \to p$ as $t \to \infty$ then ensures that the terms within the braces above tend to $\frac{p}{1-p}$ by continuity.  $\qed$

\section{Conclusion}\label{section:conclusion}
The most obvious open problem is to resolve Question A for other $\mathcal{P}$, of which increasing properties are especially natural (e.g., containment of fixed subgraphs, high chromatic number, non-planarity).  The most famous such problem is to find a bound $\mu(\text{contains $P_3$}, n) < 0.4999$, but a result of the form $\mu(\text{contains $B$}, n) < 0.4999$ for \textit{any} particular bipartite $B$ would likely be great progress.

Closer to the results of this paper, we suggest Question A in the case of Hamiltonicity.  By Theorem \ref{connected} (and from the canonical lower bound of fixing a particular Hamiltonian cycle), we have
\[
2^{-n} \leq \mu(\text{Hamiltonian}, n) \leq \mu(\text{connected}, n) = 2^{-n+1}.
\]
Here, we conjecture that the lower bound is tight, but it is unclear how to improve either bound.

It would also be interesting to study $\mu(\text{contains $P_k$}, n)$, the path with $k$ edges.  For $k=n-1$, Theorem \ref{connected} shows this is $2^{-n+1}$, and we have already discussed the case $k=3$.  A natural next step might be to study $P_{n-t}$ for other fixed $t$ and ask when $\mu(\text{contains $P_k$}, n) = 2^{-k}$.  We also suggest the properties ``has at most $k$ connected components" and also ``has minimum degree at least $k$."

Another direction of research could be to ask for quantitative improvements in Theorem \ref{unions}.  That is, for a fixed graph $\Gamma$, let $a_t = \lim_{n \to \infty} \mu(\text{contains $\Gamma_t$}, n)$.  A tensoring construction similar to that of Section \ref{section:winkler} shows $a_{k+l} \geq a_{k} a_{l}$, so by Fekete's lemma, $\lim_{t \to \infty} a_t ^{1/t}$ exists.  Theorem \ref{unions} proves this limit is at least $a_1 / (1-a_1)$, and it is natural to ask whether or not this is optimal.  However, considering the weighted version of Katona's $t$-intersection theorem (see \cite{Tintersecting} for a solution and discussion, and see \cite{winkler} for an elegant related coupling argument), we believe that any improvement to our lower bound would require a rather different approach.

\bibliographystyle{siam}
\bibliography{../mybib}

\bigskip
\myAuthor{Aaron Berger}{bergera@mit.edu}{Massachusetts Institute of Technology}

\medskip
\myAuthor{Ross Berkowitz}{ross.berkowitz@yale.edu}{Yale University}

\medskip
\myAuthor{Pat Devlin}{patrick.devlin@yale.edu}{Yale University}

\medskip
\myAuthor{Michael Doppelt}{michael.doppelt@yale.edu}{Yale University}

\medskip
\myAuthor{Sonali Durham}{sonali.durham@yale.edu}{Yale University}

\medskip
\myAuthor{Tessa Murthy}{tessa.murthy@yale.edu}{Yale University}

\medskip
\myAuthor{Harish Vemuri}{harish.vemuri@yale.edu}{Yale University}

\end{document}